\documentclass[11pt,leqno]{article}

\usepackage{amsmath,amsfonts,amscd,amssymb,theorem}

\long\def\comment#1\endcomment{}

\comment
\endcomment

\makeatletter
\begingroup
\gdef\th@dotted{\normalfont\itshape
  \def\@begintheorem##1##2{%
        \item[\hskip\labelsep \theorem@headerfont ##1\ ##2.]}%
\def\@opargbegintheorem##1##2##3{%
   \item[\hskip\labelsep \theorem@headerfont ##1\ ##2\ (##3).]}}
\endgroup
\makeatother

\theoremstyle{dotted}

\newtheorem{theorem}{Theorem}[section]
\newtheorem{lemma}[theorem]{Lemma}

\newtheorem{prop}[theorem]{Proposition}

\makeatletter
\begingroup
\gdef\th@upshape{\normalfont
  \def\@begintheorem##1##2{%
        \item[\hskip\labelsep \theorem@headerfont ##1\ ##2.]}%
\def\@opargbegintheorem##1##2##3{%
   \item[\hskip\labelsep \theorem@headerfont ##1\ ##2\ (##3).]}}
\endgroup
\makeatother

\theoremstyle{upshape}

\newtheorem{defn}[theorem]{Definition}

\makeatletter
\renewcommand{\subsection}{\@startsection{subsection}{2}{0pt}{-3ex
plus -1ex minus -0.2ex}{-2mm plus -0pt minus
-2pt}{\normalfont\bfseries}} \makeatother

\makeatletter
\@addtoreset{equation}{section}
\makeatother

\newcommand{\cntrct}                
{\hspace{2pt}\raisebox{1pt}{\text{$\lrcorner$}}\hspace{2pt}}

\newcommand{\proof}[1][Proof.]{\smallskip\noindent{\em #1}}
\def\endproof{\hfill\ensuremath{\square}\par\medskip}

\def\eqref#1{\thetag{\ref{#1}}}

\let\latexref=\ref
\def\ref#1{{\normalfont{\latexref{#1}}}}

\newcommand{\wt}{\widetilde}
\newcommand{\wh}{\widehat}

\newcommand{\idot}{{\:\raisebox{1pt}{\text{\circle*{1.5}}}}}
\newcommand{\hdot}{{\:\raisebox{3pt}{\text{\circle*{1.5}}}}}
\newcommand{\Z}{{\mathbb Z}}

\newcommand{\Q}{\overline{Q}}

\newcommand{\eps}{\varepsilon}
\renewcommand{\phi}{\varphi}

\newcommand{\calo}{{\mathcal O}}

\newcommand{\vH}{\check{H}}

\def\dlim_#1{{\displaystyle\lim_{#1}}^\hdot}

\newcommand{\Hom}{\operatorname{Hom}}
\newcommand{\Ext}{\operatorname{Ext}}

\newcommand{\Tor}{\operatorname{Tor}}
\newcommand{\Coker}{\operatorname{Coker}}
\newcommand{\Ker}{\operatorname{Ker}}
\renewcommand{\Im}{\operatorname{Im}}

\newcommand{\Fun}{\operatorname{Fun}}

\newcommand{\id}{\operatorname{\sf id}}
\newcommand{\Id}{\operatorname{\sf Id}}

\newcommand{\tw}{ {(1)} }

\newcommand{\A}{{\cal A}}
\newcommand{\V}{{\cal V}}
\newcommand{\D}{{\cal D}}

\newcommand{\C}{{\mathbb C}}

\newcommand{\Hh}{{\mathbb H}}
\newcommand{\hh}{{\mathcal H}}

\newcommand{\bimod}{\operatorname{\!-\sf bimod}}
\renewcommand{\mod}{{\text{-\sf mod}}}
\newcommand{\fg}{{\text{\sf\small fg}}}

\newcommand{\Vect}{\operatorname{\!-Vect}}

\newcommand{\Span}{\operatorname{\sf Span}}

\newcommand{\Add}{\operatorname{\sf Add}}
\newcommand{\St}{\operatorname{\sf St}}

\newcommand{\Spec}{\operatorname{Spec}}

\newcommand{\cchar}{\operatorname{\sf char}}

\title{Cartier isomorphism and Hodge Theory in the non-commutative case}

\author{D. Kaledin\thanks{Partially supported by CRDF grant RUM1-2694.}}

\begin{document}

\maketitle

\tableofcontents

\section{Introduction}

One of the standard ways to compute the cohomology groups of a
smooth complex manifold $X$ is by means of the de Rham theory: the
de Rham cohomology groups
\begin{equation}\label{derham.def}
H_{DR}^\hdot(X) = \Hh^\hdot(X,\Omega_{DR}^\hdot)
\end{equation}
are by definition the hypercohomology groups of $X$ with
coefficients in the (holomorphic) de Rham complex
$\Omega_{DR}^\hdot$, and since, by Poincar\'e Lemma,
$\Omega_{DR}^\hdot$ is a resolution of the constant sheaf $\C$, we
have $H_{DR}^\hdot(X) \cong H^\hdot(X,\C)$. If $X$ is in fact
algebraic, then $\Omega_{DR}^\hdot$ can also be defined
algebraically, so that the right-hand side in \eqref{derham.def} can
be understood in two ways: either as the hypercohomology of an
analytic space, or as the hypercohomology of a scheme equipped with
the Zariski topology. One can show that the resulting groups
$H^\hdot_{DR}(X)$ are the same (for compact $X$, this is just the
GAGA principle; in the non-compact case this is a difficult but true
fact established by Grothendieck \cite{Gr}).

Of course, an algebraic version of the Poincar\'e Lemma is false,
since the Zariski topology is not fine enough -- no matter how small
a Zariski neighbohood of a point one takes, it usually has
non-trivial higher de Rham cohomology. However, the Lemma survives
on the formal level: the completion $\wh{\Omega^\hdot_{DR}}$ of the
de Rham complex near a closed point $x \in X$ is quasiisomorphic to
$\C$ placed in degree $0$.

\smallskip

Assume now that our $X$ is a smooth algebraic variety over a perfect
field $k$ of characteristic $p > 0$. Does the de Rham cohomology
still make sense?

\smallskip

The de Rham complex $\Omega_{DR}^\hdot$ itself is well-defined:
$\Omega^1$ is just the sheaf of K\"ahler differentials, which makes
sense in any characteristic and comes equipped with the universal
derivation $d:\calo_X \to \Omega^1$, and $\Omega_{DR}^\hdot$ is its
exterior algebra, which is also well-defined in characteristic $p$.
However, the Poincar\'e Lemma breaks down completely -- the homology
of the de Rham complex remains large even after taking completion at
a closed point.

In degree $0$, this is actually very easy to see: for any local
function $f$ on $X$, we have $df^p = pf^{p-1}df = 0$, so that all
the $p$-th powers of functions are closed with respect to the de
Rham differential. Since we are in characteristic $p$, these powers
form a subsheaf of algebras in $\calo_X$ which we denote by
$\calo_X^p \subset \calo_X$. This is a large subsheaf. In fact, if
we denote by $X^\tw$ the scheme $X$ with $\calo_X^p$ as the
structure sheaf, then $X \cong X^\tw$ as abstract schemes, with the
isomorphism given by the Frobenius map $f \mapsto f^p$. Fifty years
ago P. Cartier proved that in fact {\em all} the functions in
$\calo_X$ closed with respect to the de Rham differential are
contained in $\calo_X^p$, and moreover, one has a similar
description in higher degrees: there exist natural isomorphisms
\begin{equation}\label{car.usu}
C:\hh^\hdot_{DR} \cong \Omega^\hdot_{X^\tw},
\end{equation}
where on the left we have the homology sheaves of the de Rham
complex, and on the right we have the sheaves of differential forms
on the scheme $X^\tw$. These isomorphisms are known as {\em Cartier
isomorphisms}.

\smallskip

The Cartier isomorphism has many applications, but one of the most
unexpected has been discovered in 1987 by P. Deligne and L. Illusie:
one can use the Cartier isomorphism to give a purely algebraic proof
of the following purely algebraic statement, which is normally
proved by the highly trancendental Hodge Theory.

\begin{theorem}[\cite{DI}]\label{di.thm}
Assum given a smooth proper variety $X$ over a field $K$ of
characteristic $0$. Then the Hodge-to-de Rham spectal sequence
$$
H^\hdot(X,\Omega^\hdot) \Rightarrow H_{DR}^\hdot(X)
$$
associated to the stupid filtration on the de Rham complex
$\Omega^\hdot$ degenerates at the first term.
\end{theorem}

The proof of Deligne and Illusie was very strange, because it worked
by reduction to positive characteristic, where the statement is not
true for a general $X$. What they proved is that if one imposes two
additional conditions on $X$, then the Cartier isomorphisms can be
combined together into a quasiisomorphism
\begin{equation}\label{car.eq}
\Omega_{DR}^\hdot \cong \bigoplus_i\hh_{DR}^i[-i] \cong
\bigoplus_i\Omega_{X^\tw}^i[-i]
\end{equation}
in the derived category of coherent sheaves on $X^\tw$.  The
degeneration follows from this immediately for dimension
reasons. The additional conditions are:
\begin{enumerate}
\item $X$ can be lifted to a smooth scheme over $W_2(k)$, the ring
  of second Witt vectors of the perfect field $k$ (e.g. if
  $k=\Z/p\Z$, $X$ has to be liftable to $\Z/p^2\Z$), and
\item we have $p > \dim X$.
\end{enumerate}
To deduce Theorem~\ref{di.thm}, one finds by the standard argument a
proper smooth model $X_R$ of $X$ defined over a finitely generated
subring $R \subset K$, one localizes $R$ so that it is unramified
over $\Z$ and all its residue fields have characteristic greater
than $\dim X$, and one deduces that all the special fibers of
$X_R/R$ satisfy the assumptions above; hence the differentials in
the Hodge-to-de Rham spectral sequence vanish at all closed points
of $\Spec R$, which means they are identically $0$ by Nakayama.

\bigskip

The goal of these lectures is to present in a down-to-earth way the
results of two recent papers \cite{K1}, \cite{K2}, where the story
summarized above has been largely transferred to the setting of {\em
non-commutative geometry}.

\bigskip

To explain what I mean by this, let us first recall that a
non-commu\-ta\-ti\-ve version of differential forms has been known
for quite some time now. Namely, assume given an associative unital
algebra $A$ over a field $k$, and an $A$-bimodule $M$. Then its {\em
Hochschild homology} $HH_\idot(A,M)$ of $A$ with coefficients in $M$
is defined as
\begin{equation}\label{hh.eq}
HH_\idot(A) = \Tor^\hdot_{A^{opp} \otimes A}(A,M),
\end{equation}
where $A^{opp} \otimes A$ is the tensor product of $A$ and the
opposite algebra $A^{opp}$, and the $A$-bimodule $M$ is treated as a
left module over $A^{opp} \otimes A$. Hochschild homology
$HH_\idot(A)$ is the Hochschild homology of $A$ with coefficients in
itself.

Assume for a moment that $A$ is in fact commutative, and $\Spec A$
is a smooth algebraic variety over $K$. Then it has been proved back
in 1962 in the paper \cite{HKR} that we have canonical isomorphisms
$HH_i(A) \cong \Omega^i(A/k)$ for any $i \geq 0$. Thus for a general
$A$, one can treat Hochschild homology classes as a replacement for
differential forms.

Moreover, in the early 1980-es it has been discovered by A.~Connes
\cite{C}, J.-L.~Loday and D.~Quillen \cite{LQ}, and B.~Feigin and
B.~Tsygan \cite{FT1}, that the de Rham differential also makes sense
in the general non-commutative setting. Namely, these authors
introduce a new invariant of an associative algebra $A$ called {\em
cyclic homology}; cyclic homology, denoted $HC_\idot(A)$, is related
to the Hochschild homology $HH_\idot(A)$ by a spectral sequence
\begin{equation}\label{h.dr.eq}
HH_\idot(A)[u^{-1}] \Rightarrow HC_\idot(A),
\end{equation}
which in the smooth commutative case reduces to the Hodge-to-de Rham
spectral sequence (here $u$ is a formal parameter of cohomological
degree $2$, and $HH_\idot(A)[u^{-1}]$ is shorthand for ``polynomials in
$u^{-1}$ with coefficients in $HH_\idot(A)$'').

It has been conjectured for some time now that the spectral sequence
\eqref{h.dr.eq}, or a version of it, degenerates under appropriate
assumptions on $A$ (which imitate the assumptions of
Theorem~\ref{di.thm}). Following \cite{K2}, we will attack this
conjecture by the method of Deligne and Illusie. To do this, we will
introduce a certain non-commutative version of the Cartier
isomorphism, or rather, of the ``globalized'' isomorphism
\eqref{car.eq} (in the process of doing it, we will need to
introduce some conditions on $A$ which precisely generalize the
conditions \thetag{i}, \thetag{ii} above). Then we prove a version
of the degeneration conjecture as stated by M. Kontsevich and
Ya. Soibelman in \cite{KS} (we will have to impose an additional
technical assumption which, fortunately, is not very drastic).

\medskip

The paper is organized as follows. In Section~\ref{cycl.sec} we
recall the definition of the cyclic homology and some versions of it
needed for the Cartier isomorphism (most of this material is quite
standard, the reader can find good expositions in \cite{L} or
\cite{FT}). One technical result needed in the main part of the
paper has been separated into Section~\ref{vani.sec}. In
Section~\ref{qfr.sec}, we construct the Cartier isomorphism for an
algebra $A$ equipped with some additional piece of data which we
call the {\em quasi-frobenius map}. It exists only for special
classes of algebras -- e.g. for free algebras, or for the group
algebra $k[G]$ of a finite group $G$ -- but the construction
illustrates nicely the general idea. In Section~\ref{add.sec}, we
show what to do in the general case. Here the conditions \thetag{i},
\thetag{ii} emerge, and in a somewhat surprising way -- as it turns
out, they essentially come from alegbraic topology, and the whole
theory has a distinct topological flavor. Finally, in
Section~\ref{deg.sec} we show how to apply our generalized Cartier
isomorphism to the Hodge-to-de Rham degeneration. The exposition in
Sections~\ref{cycl.sec}-\ref{qfr.sec} is largely self-contained. In
the rest of the paper, we switch to a more descriptive style, with
no proofs, and not many precise statements; this part of the paper
should be treated as a companion to \cite{K2}.

\subsection*{Acknowledgements.}
This paper is a write-up (actually quite an enlarged write-up) of
two lectures given in Goettingen in August 2006, at a summer school
organized by Yu.~Tschinkel and funded by the Clay Institute. I am
very grateful to all concerned for making it happen, and for giving
me an opportunity to present my results. In addition, I would like
to mention that a large part of the present paper is written in
overview style, and many, if not most of the things overviewed are
certainly {\em not} my results. This especially concerns
Section~\ref{cycl.sec}, on one hand, and Section~\ref{deg.sec}, on
the other hand. Given the chosen style, it is difficult to provide
exact attributions; however, I should at least mention that I've
learned much of this material from A. Beilinson, A. Bondal,
M. Kontsevich, B. To\"en and B. Tsygan.

\section{Cyclic homology package}\label{cycl.sec}

\subsection{Basic definitions.}

The fastest way and most down-to-earth to define cyclic homology is
by means of an explicit complex. Namely, assume given an associative
unital algebra $A$ over a field $k$. To compute its Hochschild
homology with coeffients in some bimodule $M$, one has to find a
flat resolution of $M$. One such is the {\em bar resolution} -- it
is rather inconvenient in practical computations, but it is
completely canonical, and it exists without any assumptions on $A$
and $M$. The terms of this resolution are of the form $A^{\otimes n}
\otimes M$, $n \geq 0$, and the differential $b':A^{\otimes n+1}
\otimes M \to A^{\otimes n} \otimes M$ is given by
\begin{equation}\label{bpr.eq}
b' = \sum_{0 \leq i \leq n} (-1)^i\id^{\otimes i} \otimes m \otimes
\id^{\otimes n-i},
\end{equation}
where $m:A \otimes A \to A$, $m:A \otimes M \to M$ are the
multiplication maps. Substituting this resolution into \eqref{hh.eq}
gives a complex which computes $HH_\idot(A,M)$; its terms are also
$A^{\otimes i} \otimes M$, but the differential is given by
\begin{equation}\label{b.eq}
b = b' + (-1)^{n+1}t,
\end{equation}
with the correction term $t$ being equal to
$$
t(a_0 \otimes \dots \otimes a_{n+1} \otimes m) = a_1 \otimes \dots
\otimes a_{n+1} \otimes ma_0
$$
for any $a_0,\dots,a_{n+1} \in A$, $m \in M$. Geometrically, one can
think of the components $a_0,\dots,a_{n-1},m$ of some tensor in
$A^{\otimes n} \otimes M$ as having been placed at $n+1$ points on
the unit interval $[0,1]$, including the egde points $0,1 \in
[0,1]$; then each of the terms in the differential $b'$ corresponds
to contracting an interval between two neighboring points and
multiplying the components sitting at its endpoints. To visualize
the differential $b$ in a similar way, one has to take $n+1$ points
placed on the unit circle $S^1$ instead of the unit interval,
including the point $1 \in S^1$, where we put the component $m$.

In the case $M = A$, the terms in the bar complex are just
$A^{\otimes n+1}$, $n \geq 0$, and they acquire an additional
symmetry: we let $\tau:A^{\otimes n+1} \to A^{\otimes n+1}$ to be
the cyclic permutation multiplied by $(-1)^n$. Note that in spite of
the sign change, we have $\tau^{n+1} = \id$, so that it generates an
action of the cyclic group $\Z/(n+1)\Z$ on every $A^{\otimes
n+1}$. The fundamental fact here is the following.

\begin{lemma}[{\cite{FT},\cite{L}}]\label{cycl.lemma}
For any $n$, we have
\begin{align*}
(\id - \tau) \circ b' &= -b \circ (\id - \tau),\\
(\id + \tau + \dots + \tau^{n-1}) \circ b &= -b' \circ (\id +
\tau + \dots + \tau^n)
\end{align*}
as maps from $A^{\otimes n+1}$ to $A^{\otimes n}$.
\end{lemma}

\proof{} Denote $m_i = \id^i \otimes m \otimes \id^{n-i}:A^{\otimes
n+1} \to A^{\otimes n}$, $0 \leq i \leq n-1$, so that $b' = m_0 -
m_1 + \dots + (-1)^{n-1}m_{n-1}$, and let $m_n = t = (-1)^n(b -
b')$. Then we obviously have
$$
m_{i+1} \circ \tau = \tau \circ m_i
$$
for $0 \leq i \leq n-1$, and $m_0 \circ \tau = (-1)^nm_n$. Formally
applying these
identities, we conclude that
\begin{equation}\label{hc.bim.1}
\begin{aligned}
\sum_{0 \leq i \leq n}(-1)^im_i \circ (\id - \tau) 
&= \sum_{0 \leq i \leq n}(-1)^im_i - m_0 - \sum_{1 \leq i \leq
n}(-1)^i\tau \circ m_{i-1} \\
&= -(\id - \tau) \circ \sum_{0 \leq i \leq n-1}(-1)^im_i,
\end{aligned}
\end{equation}
\begin{multline}\label{hc.bim.2}
b' \circ (\id + \tau + \dots + \tau^n) = \sum_{0 \leq i \leq n-1}
\sum_{0 \leq j \leq n} (-1)^im_i \circ \tau^j \\
\begin{aligned}
&= \sum_{0 \leq j \leq i \leq n-1}(-1)^i\tau^j
\circ m_{i-j} + \sum_{1 \leq i \leq j \leq n}
(-1)^{i+n}\tau^{j-1} \circ m_{n+i-j} \\
&= -(\id + \tau + \dots + \tau^{n-1}) \circ b,
\end{aligned}
\end{multline}
which proves the claim.
\endproof

As a corollary, the following diagram is in fact a bicomplex.
\begin{equation}\label{hc.eq}
\begin{CD}
\dots @>>> A @>{\id}>> A @>{0}>> A\\
@. @AA{b}A @AA{b'}A @AA{b}A \\
\dots @>>> A \otimes A @>{\id + \tau}>> A \otimes A 
@>{\id - \tau}>> A \otimes A\\
@. @AA{b}A @AA{b'}A @AA{b}A\\
\dots @. \dots @. \dots @. \dots\\
@. @AA{b}A @AA{b'}A @AA{b}A \\
\dots @>>> A^{\otimes n} @>{\id + \tau + \dots + \tau^{n-1}}>> 
A^{\otimes n} @>{\id - \tau}>> A^{\otimes n}\\
@. @AA{b}A @AA{b'}A @AA{b}A
\end{CD}
\end{equation}
Here it is understood that the whole thing extends indefinitely to
the left, all the even-numbered columns are the same, all
odd-numbered columns are the same, and the bicomplex is invariant
with respect to the horizontal shift by $2$ columns. The total
homology of this bicomplex is called the {\em cyclic homology} of
the algebra $A$, and denoted by $HC_\idot(A)$.

We see right away that the first, the third, and so on column when
counting from the right is the bar complex which computes
$HH_\idot(A)$, and the second, the fourth, and so on column is
acyclic (the top term is $A$, and the rest is the bar resolution for
$A$). Thus the spectral sequence for this bicomplex has the form
given in \eqref{h.dr.eq} (modulo obvious renumbering). On the other
hand, the rows of the bicomplex are just the standard $2$-periodic
complexes which compute the cyclic group homology
$H_\idot(\Z/n\Z,A^{\otimes n})$ (with respect to the $\Z/n\Z$-action
on $A^{\otimes n}$ given by $\tau$).

Shifting \eqref{hc.eq} to the right by $2$ columns gives the {\em
periodicity map} $u:HC_{\idot+2}(A) \to HC_\idot(A)$, which fits
into an exact triangle
\begin{equation}\label{connes.ex}
\begin{CD}
HH_{\idot+2} @>>> HC_{\idot+2}(A) @>>> HC_\idot(A) @>>>,
\end{CD}
\end{equation}
known as the {\em Connes' exact sequence}. One can also invert the
periodicity map -- in other words, extend the bicomplex
\eqref{hc.eq} not only to the left, but also to the right. This
gives the {\em periodic cyclic homology} $HP_\idot(A)$. Since the
bicomplex for $HP_\idot(A)$ is infinite in both directions, there is
a choice involved in taking the total complex: we can take either
the product, or the sum of the terms. We take the product. In
characteristic $0$, the sum is actually acyclic (because so is every
row).

If $A$ is commutative, $X = \Spec(A)$ is smooth, and $\cchar k$ is
either $0$ or greater than $\dim X$, then the only non-trivial
differential in the Hodge-to-de Rham spectral sequence
\eqref{h.dr.eq} is the first one, and it is the de Rham
differential. Consequently, we have $HP_\idot(A) =
H^\hdot_{DR}(X)((u))$ (where as before, $u$ is a formal variable of
cohomological degree $2$).

\subsection{The $p$-cyclic complex.}\label{p.cycl.subs}

All of the above is completely standard; however, we will also need
to use another way to compute $HC_\idot(A)$, which is less
standard. Namely, fix an integer $p \geq 2$, and consider the
algebra $A^{\otimes p}$. Let $\sigma:A^{\otimes p} \to A^{\otimes
p}$ be the cyclic permutation, and let $A^{\otimes p}_\sigma$ be the
diagonal $A^{\otimes p}$-bimodule with the bimodule structure
twisted by $\sigma$ --namely, we let
$$
a \cdot b \cdot c = ab\sigma(c)
$$
for any $a,b,c \in A^{\otimes p}$.

\begin{lemma}\label{hh.p}
We have $HH_\idot(A^{\otimes p},A^{\otimes p}_\sigma) \cong
HH_\idot(A)$.
\end{lemma}

\proof{} Induction on $p$. We may compute the tensor product in
\eqref{hh.eq} over each of the factors $A$ in $A^{\otimes p}$ in
turn; this shows that
$$
HH_\idot(A^{\otimes p},A^{\otimes p}_\sigma) \cong
\Tor^\hdot_{(A^{\otimes (p-1)})^{opp} \otimes A^{\otimes
    (p-1)}}\left(A^{\otimes (p-1)},\Tor^\hdot_{A^{opp} \otimes
  A}(A,A^{\otimes p}_\sigma)\right),
$$
and one check easily that as long as $p \geq 2$, so that $A^{\otimes
p}_\sigma$ is flat over $A^{opp} \otimes A$, $\Tor^i_{A^{opp}
\otimes A}(A,A^{\otimes p}_\sigma)$ is naturally isomorphic to
$A^{\otimes (p-1)}_\sigma$ if $i=0$, and trivial if $i \geq 1$.
\endproof

By virtue of this Lemma, we can use the bar complex for the algebra
$A^{\otimes p}$ to compute $HH_\idot(A)$. The resulting complex has
terms $A^{\otimes pn}$, $n \geq 0$. The differential
$b_p':A^{\otimes p(n+2)} \to A^{\otimes p(n+1)}$ is given by
essentially the same formula as \eqref{bpr.eq}:
$$
b'_p = \sum_{0 \leq i \leq n} (-1)^im_i^p = \sum_{0 \leq i \leq n}
(-1)^i\id^{\otimes pi} \otimes m^{\otimes p} \otimes \id^{\otimes
p(n-i)},
$$
where we decompose $A^{\otimes p(n+1)} = \left(A^{\otimes
p}\right)^{\otimes (n+1)}$. The correcting term $t_p = m^p_{n+1}$ in
$\eqref{b.eq}$ is given by $m_0 \circ \tau$ (where, as before,
$\tau:A^{\otimes p(n+2)}$ is the cyclic permutation of order
$p(n+2)$ twisted by a sign). Geometrically, the component $m^p_i$ of
the Hochschild differential $b_p$ correspond to contracting
simultaneously the $i$-th, the $(i+p)$-th, the $(i+2p)$-th, and so
on interval in the unit circle divived into $p(n+2)$ intervals by
$p(n+2)$ points. On the level of bar complexes, the comparison
isomorphism $HH_\idot(A^{\otimes p},A^{\otimes p}_\sigma) \cong
HH_\idot(A)$ of Lemma~\ref{hh.p} is represented by the map 
\begin{equation}\label{M.eq}
M = m \circ (id \otimes m) \circ (\id^{\otimes 2} \otimes m) \circ
\dots \circ (\id^{\otimes pn-2} \otimes m):A^{\otimes pn} \to
A^{\otimes n};
\end{equation}
explicitly, we have
\begin{multline*}
M(a_{1,1} \otimes a_{2,1} \otimes \dots \otimes a_{n,1} \otimes
a_{1,2} \otimes a_{2,2} \otimes \dots \otimes a_{n,2} \otimes \dots
\otimes a_{1,p} \otimes a_{2,p} \otimes \dots \otimes a_{n,p}) \\
=
a_{1,1} \otimes a_{2,1} \otimes \dots \otimes a_{n-1,1} \otimes
\left(a_{n,1} \cdot \prod_{2 \leq j \leq p}\prod_{1 \leq i \leq
  n}a_{i,j}\right)
\end{multline*}
for any $a_{1,1} \otimes a_{2,1} \otimes \dots \otimes a_{n,1}
\otimes a_{1,2} \otimes a_{2,2} \otimes \dots \otimes a_{n,2}
\otimes \dots \otimes a_{1,p} \otimes a_{2,p} \otimes \dots \otimes
a_{n,p} \in A^{\otimes pn}$ -- in other words, $M:A^{\otimes pn} \to
A^{\otimes n}$ leaves the first $n-1$ terms in the tensor product
intact and multiplies the remaining $pn-n+1$ terms. We leave it to
an interested reader to check explicitly that $M \circ b_p = b \circ
M$.

\begin{lemma}\label{p.cycl.lemma}
For any $n$, we have
\begin{align*}
(\id - \tau) \circ b'_p &= -b_p \circ (\id - \tau),\\
(\id + \tau + \dots + \tau^{pn-1}) \circ b_p &= -b'_p \circ (\id +
\tau + \dots + \tau^{p(n+1)-1})
\end{align*}
as maps from $A^{\otimes p(n+1)}$ to $A^{\otimes pn}$.
\end{lemma}

\proof{} One immediately checks that, as in the proof of
Lemma~\ref{cycl.lemma}, we have
$$
m^p_{i+1} \circ \tau = - \tau \circ m^p_i
$$
for $0 \leq i \leq n$, and we also have $m^p_0 \circ \tau =
m^p_{n+1}$. Then the first equality follows from \eqref{hc.bim.1},
and \eqref{hc.bim.2} gives
$$
(\id + \tau + \dots + \tau^{n-1}) \circ b_p = -b'_p \circ (\id +
\tau + \dots + \tau^n)
$$
(note that the proof of these two equalities does {\em not} use the
fact that $\tau^{n+1}=\id$ on $A^{\otimes(n+1)}$). To deduce the
second equality of the Lemma, it suffices to notice that
$$
\id + \tau + \dots + \tau^{p(n+1)-1} = (\id + \tau + \dots + \tau^n)
\circ (\id + \sigma + \dots + \sigma^{p-1}),
$$
and $\sigma$ commutes with all the maps $m^p_i$.
\endproof

Using Lemma~\ref{p.cycl.lemma}, we can construct a version of the
bicomplex \eqref{hc.eq} for $p > 1$:
\begin{equation}\label{hc.eq.p}
\begin{CD}
\dots @>>> A^{\otimes p} @>{\id + \tau + \dots + \tau^{p-1}}>>
A^{\otimes p} @>{\id - \tau}>> A^{\otimes p}\\ 
@. @AA{b_p}A @AA{b_p'}A @AA{b_p}A \\
\dots @>>> A^{\otimes 2p} @>{\id + \dots + \tau^{2p-1}}>> A^{\otimes 2p}
@>{\id - \tau}>> A^{\otimes 2p}\\
@. @AA{b_p}A @AA{b_p'}A @AA{b_p}A\\
\dots @. \dots @. \dots @. \dots\\
@. @AA{b_p}A @AA{b_p'}A @AA{b_p}A \\
\dots @>>> A^{\otimes pn} @>{\id + \tau + \dots + \tau^{pn-1}}>> 
A^{\otimes pn} @>{\id - \tau}>> A^{\otimes pn}\\
@. @AA{b_p}A @AA{b_p'}A @AA{b_p}A
\end{CD}
\end{equation}
By abuse of notation, we denote the homology of the total complex of
this bicomplex by $HC_\idot(A^{\otimes p},A^{\otimes p}_\sigma)$.
(This is really abusive, since in general one {\em cannot} define
cyclic homology with coefficients in a bimodule -- unless the
bimodule is equipped with additional structure, as e.g.\ in
\cite{K3}, which lies beyond the scope of this paper.) As for the
usual cyclic complex, we have the periodicity map, the Connes' exact
sequence, and we can form the periodic cyclic homology
$HP_\idot(A^{\otimes p},A^{\otimes p}_\sigma)$.

\subsection{Small categories.}

Unfortunately, this is as far as the down-to-earth approach takes
us. While it is true that the isomorphism $HH_\idot(A^{\otimes
p},A^{\otimes p}_\sigma) \cong HH_\idot(A)$ given in
Lemma~\ref{hh.p} can be extended to an isomorphism
$$
HC_\idot(A^{\otimes p},A^{\otimes p}_\sigma) \cong HC_\idot(A),
$$
it is not possible to realize this extended isomorphism by an
explicit map of bicomplexes. Indeed, already in degree $0$ the
comparison map $M$ of \eqref{M.eq} which realized the isomorphism
$$
HH_0(A^{\otimes p},A^{\otimes p}_\sigma) \to HH_0(A)
$$
on the level of bar complexes is given by the multiplication map
$A^{\otimes p} \to A$, and to define this multiplication map, one
has to break the cyclic symmetry of the product $A^{\otimes p}$. The
best one can obtain is a map between total complexes computing
$HC_\idot(A^{\otimes p},A^{\otimes p}_\sigma)$ and $HC_\idot(A)$
which preserves the filtration, but not the second grading; when one
tries to write the map down explicitly, the combinatorics quickly
gets completely out of control.

For this reason, in \cite{K1} and \cite{K2} one follows \cite{C} and
uses a more advanced approach to cyclic homology which is based on
the technique of {\em homology of small categories} (see
e.g. \cite[Section 6]{L}). Namely, for any small category $\Gamma$
and any base field $k$, the category $\Fun(\Gamma,k)$ of functors
from $\Gamma$ to $k$-vector spaces is an abelian category, and the
direct limit functor $\lim_{\overset{\to}{\Gamma}}$ is
right-exact. Its derived functors are called {\em homology functors}
of the category $\Gamma$ and denoted by $H_\idot(\Gamma,E)$ for any
$E \in \Fun(\Gamma,k)$. For instance, if $\Gamma$ is a groupoid with
one object with automorphism group $G$, then $\Fun(\Gamma,k)$ is the
category of $k$-representations of the group $G$; the homology
$H_\idot(\Gamma,-)$ is then tautologically the same as the group
homology $H_\idot(G,-)$. Another example is the category
$\Delta^{opp}$, the opposite to the category $\Delta$ of finite
non-empty totally-ordered sets. It is not difficult to check that
for any simplicial $k$-vector $E \in \Fun(\Delta^{opp},k)$, the
homology $H_\idot(\Delta^{opp},E)$ can be computed by the standard
chain complex of $E$.

For applications to cyclic homology, one introduces special small
categories $\Lambda_\infty$ and $\Lambda_p$, $p \geq 1$. The objects
in the category $\Lambda_\infty$ are numbered by positive integers
and denoted $[n]$, $n \geq 1$. For any $[n],[m] \in \Lambda_\infty$,
the set of maps $\Lambda_\infty([n],[m])$ is the set of all maps
$f:\Z \to \Z$ such that
\begin{equation}\label{lambda.defn.eq}
\begin{aligned}
f(a) \leq f(b) \quad \text{ whenever }a \leq b,\qquad
f(a+n) = f(a)+m,
\end{aligned}
\end{equation}
for any $a,b \in \Z$. For any $[n] \in \Lambda_\infty$, denote by
$\sigma:[n] \to [n]$ the endomorphism given by $f(a)=a+n$. Then
$\sigma$ commutes with all maps in $\Lambda_\infty$. The category
$\Lambda_p$ has the same objects as $\Lambda_\infty$, and the set of
maps is
$$
\Lambda_p([n],[m]) = \Lambda_\infty([n],[m])/\sigma^p
$$
for any $[n],[m] \in \Lambda_p$. The category $\Lambda_1$ is denoted
simply by $\Lambda$; this is the original cyclic category
introduced by A. Connes in \cite{C}. By definition, we have
projections $\Lambda_\infty \to \Lambda_p$ and $\pi:\Lambda_p \to
\Lambda$.

If we only consider those maps in \eqref{lambda.defn.eq} which send
$0 \in \Z$ to $0$, then the resulting subcategory in
$\Lambda_\infty$ is equivalent to $\Delta^{opp}$. This gives a
canonical embedding $j:\Delta^{opp} \to \Lambda_\infty$, and
consequently, embeddings $j:\Delta^{opp} \to \Lambda_p$.

The category $\Lambda_p$ conveniently encodes the maps $m^p_i$ and
$\tau$ between various tensor powers $A^{\otimes pn}$ used in the
complex \eqref{hc.eq.p}: $m^p_i$ corresponds to the map $f \in
\Lambda_p([n+1],[n])$ given by
$$
f(a(n+1)+b) =
\begin{cases}
an+b, &\quad b \leq i,\\
an+b-1, &\quad b > i,
\end{cases}
$$
where $0 \leq b \leq n$, and $\tau$ is the map $a \mapsto a+1$,
twisted by the sign (alternatively, one can say that $m_i^p$ are
obtained from face maps in $\Delta^{opp}$ under the embedding
$\Delta^{opp} \subset \Lambda_p$). The relations between these maps
which we used in the proof of Lemma~\ref{p.cycl.lemma} are encoded
in the composition laws of the category $\Lambda_p$. Thus for any
object $E \in \Fun(\Lambda_p,k)$ -- they are called {\em $p$-cyclic
objects} -- one can form the bicomplex of the type \eqref{hc.eq.p}
(or \eqref{hc.eq}, for $p=1$):
\begin{equation}\label{hc.lambda.p}
\begin{CD}
\dots @>>> E([1]) @>{\id + \tau + \dots + \tau^{p-1}}>>
E([1]) @>{\id - \tau}>> E([1])\\ 
@. @AA{b_p}A @AA{b_p'}A @AA{b_p}A \\
\dots @>>> E([2]) @>{\id + \dots + \tau^{2p-1}}>> E([2])
@>{\id - \tau}>> E([2])\\
@. @AA{b_p}A @AA{b_p'}A @AA{b_p}A\\
\dots @. \dots @. \dots @. \dots\\
@. @AA{b_p}A @AA{b_p'}A @AA{b_p}A \\
\dots @>>> E([n]) @>{\id + \tau + \dots + \tau^{pn-1}}>> 
E([n]) @>{\id - \tau}>> E([n])\\
@. @AA{b_p}A @AA{b_p'}A @AA{b_p}A
\end{CD}
\end{equation}
Just as for the complex \eqref{hc.eq.p}, we have periodicity, the
periodic version of the complex, and the Connes' exact sequence
\eqref{connes.ex} (the role of Hochschild homology is played by the
standard chain complex of the simplicial vector space $j^*E \in
\Fun(\Delta^{opp},k)$).

\begin{lemma}\label{cycl.comp}
For any $E \in \Fun(\Lambda_p,k)$, the homology
$H_\idot(\Lambda_p,E)$ can be computed by the bicomplex
\eqref{hc.lambda.p}.
\end{lemma}

\proof{} The homology of the total complex of \eqref{hc.lambda.p} is
obviously a homological functor from $\Fun(\Lambda_p,k)$ to
$k$ (that is, short exact sequences in $\Fun(\Lambda_p,k)$
gives long exact sequences of homology). Therefore it suffices to
prove the claim for a set of projective generators of the category
$\Fun(\Lambda_p,k)$. For instance, it suffices to consider all
the representable functors $E_n$, $n \geq 1$ -- that is, the
functors given by
$$
E_n([m]) = k\left[\Lambda_p([n],[m])\right],
$$
where in the right-hand side we take the $k$-linear span. Then on
one hand, for general tautological reasons -- essentially by Yoneda
Lemma -- $H_\idot(\Lambda_p,E_n)$ is $k$ in degree $0$ and $0$ in
higher degrees. On the other hand, the action of the cyclic group
$\Z/pm\Z$ generated by $\tau \in \Lambda_p([m],[m])$ on
$\Lambda_p([n],[m])$ is obviously free, and we have
$$
\Lambda_p([n],[m])/\tau \cong \Delta^{opp}([n],[m])
$$
-- every $f:\Z \to \Z$ can be uniquely decomposed as $f = \tau^j
\circ f_0$, where $0 \leq j < pm$, and $f_0$ sends $0$ to
$0$. The rows of the complex \eqref{hc.lambda.p} compute
$$
H_\idot(\Z/pm\Z,E_n([m])) \cong k\left[\Delta^{opp}([n],[m])\right],
$$
and the first term in the corresponding spectral sequence is the
standard complex for the simplicial vector space $E_n^\Delta \in
\Fun(\Delta^{opp},k)$ represented by $[n] \in
\Delta^{opp}$. Therefore this complex computes
$H_\idot(\Delta^{opp},E_n^\Delta)$, which is again $k$.
\endproof

The complex \eqref{hc.eq} is the special case of \eqref{hc.lambda.p}
for $p=1$ and the following object $A_\# \in \Fun(\Lambda,k)$:
we set $A_\#([n]) = A^{\otimes n}$, where the factors are numbered
by elements in the set $V([n]) = \Z/n\Z$, and any $f \in
\Lambda([n],[m])$ acts by
$$
A_\#(f)\left(\bigotimes_{i \in V([n])}a_i\right)=
\bigotimes_{j \in V([m])} \prod_{i \in f^{-1}(j)}a_i,
$$
(if $f^{-1}(i)$ is empty for some $i \in V([n])$, then the
right-hand side involves a product numbered by the empty set; this
is defined to be the unity element $1 \in A$). To obtain the complex
\eqref{hc.eq.p}, we note that for any $p$, we have a functor
$i:\Lambda_p \to \Lambda$ given by $[n] \mapsto [pn]$, $f \mapsto
f$. Then \eqref{hc.lambda.p} applied to $i^*A_\# \in
\Fun(\Lambda_p,k)$ gives \eqref{hc.eq.p}. By Lemma~\ref{cycl.comp}
we have
\begin{align*}
HC_\idot(A) &\cong H_\idot(\Lambda,k),\\
HC_\idot(A^{\otimes p},A^{\otimes p}_\sigma) &\cong
H_\idot(\Lambda_p,k).
\end{align*}

\begin{lemma}[{{\cite[Lemma 1.12]{K2}}}]\label{i.lemma}
For any $E \in \Fun(\Lambda,k)$, we have a natural isomorphism
$$
H_\idot(\Lambda_p,i^*E) \cong H_\idot(\Lambda,E),
$$
which is compatible with the periodicity map and with the Connes' exact
sequence \eqref{connes.ex}.\endproof
\end{lemma}

Thus $HC_\idot(A^{\otimes p},A^{\otimes p}_\sigma) \cong
HC_\idot(A)$. The proof of this Lemma is not difficult. First of
all, a canonical comparison map $H_\idot(\Lambda_p,i^*E) \to
H_\idot(\Lambda,E)$ exists for tautological adjunction
reasons. Moreover, the periodicity homomorphism for
$H_\idot(\Lambda_p,-)$ is induced by the action of a canonical
element $u_p \in H^2(\Lambda_p,k) = \Ext^2(k,k)$, where $k$ means
the constant functor $[n] \mapsto k$ from $\Lambda_p$ to
$k$. One check explicitly that $i^*u = u_p$, so that the
comparison map is indeed compatible with periodicity, and then it
suffices to prove that the comparison map
$$
H_\idot(\Delta^{opp},i^*E) \to H_\idot(\Delta^{opp},E)
$$
is an isomorphism. When $E$ is of the form $A_\#$, this is
Lemma~\ref{hh.p}; in general, one shows that $\Fun(\Lambda,k)$
has a projective generator of the form $A_\#$. For details, we refer
the reader to \cite{K2}.

\section{One vanishing result}\label{vani.sec}

For our construction of the Cartier map, we will need one
vanishing-type result on periodic cyclic homology in prime
characteristic -- we want to claim that the periodic cyclic homology
$HP_\idot(E)$ of a $p$-cyclic object $E$ vanishes under some
assumptons on $E$.

First, consider the cyclic group $\Z/np\Z$ for some $n,p \geq 1$,
with the subgroup $\Z/p\Z \subset \Z/pn\Z$ and the quotient $\Z/n\Z
= (\Z/pn\Z)/(\Z/p\Z)$. It is well-known that for any representation
$V$ of the group $\Z/pn\Z$, we have the Hochschild-Serre spectral
sequence
$$
H_\idot(\Z/n\Z,H_\idot(\Z/p\Z,V))\Rightarrow H_\idot(\Z/pn\Z,-).
$$
To see it explicitly, one can compute the homology
$H_\idot(\Z/np\Z,V)$ by a complex which is slightly more complicated
than the standard one. Namely, write down the diagram
\begin{equation}\label{hs.eq}
\begin{CD}
@>{\id - \sigma}>> V @>{d_\sigma}>> V @>{\id -
  \sigma}>> V\\
@. @AA{\id-\tau}A @AA{\id-\tau}A @AA{\id-\tau}A\\
@>{\id - \sigma}>> V @>{d_\sigma}>> V @>{\id -
  \sigma}>> V\\
@. @AA{d_\tau}A @AA{d_\tau}A @AA{d_\tau}A\\
@>{\id - \sigma}>> V @>{d_\sigma}>> V @>{\id -
  \sigma}>> V\\
@. @AA{\id-\tau}A @AA{\id-\tau}A @AA{\id-\tau}A,
\end{CD}
\end{equation}
where $\tau$ is the generator of $\Z/pn\Z$, $\sigma=\tau^n$ is the
generator of $\Z/p\Z \subset \Z/pn\Z$, and
$d_\sigma=\id+\sigma+\dots+\sigma^p$,
$d_\tau=\id+\tau+\dots+\tau^{n-1}$. This is not quite a bicomplex
since the vertical differential squares to $\id - \sigma$, not to
$0$; to correct this, we add to the total differential the term
$\id:V \to V$ of bidegree $(-1,2)$ in every term in the columns with
odd numbers (when counting from the right). The result is a filtered
complex which computes $H_\idot(\Z/pn\Z,V)$, and the
Hochschild-Serre spectral sequence appears as the spectral sequence
of the filtered complex \eqref{hs.eq}.

One feature which is apparent in the complex \eqref{hs.eq} is that
it has two different periodicity endomorphisms: the endomorphism
which shift the diagram to the left by two columns (we will denote
it by $u$), and the endomorphism which shifts the diagram downwards
by two rows (we will denote it by $u'$).

Assume now given a field $k$ and a $p$-cyclic object $E \in
\Fun(\Lambda_p,k)$, and consider the complex
\eqref{hc.lambda.p}. Its $n$-th row is the standard periodic complex
which computes $H_\idot(\Z/pn\Z,E([n]))$, and we can replace all
these complexes by the corresponding complex \eqref{hs.eq}. By
virtue of Lemma~\ref{p.cycl.lemma}, the result is a certain filtered
bicomplex of the form
\begin{equation}\label{hc.hs.eq}
\begin{CD}
@>{\id - \sigma}>> C_\idot(E) @>{\id+\sigma+\dots+\sigma^{p-1}}>>
  C_\idot(E) @>{\id - \sigma}>> C_\idot(E)\\
@. @AA{B}A @AA{B}A @AA{B}A\\
@>{\id - \sigma}>> C'_\idot(E) @>{\id+\sigma+\dots+\sigma^{p-1}}>>
  C'_\idot(E) @>{\id - \sigma}>> C'_\idot(E)\\
@. @AA{B}A @AA{B}A @AA{B}A\\
@>{\id - \sigma}>> C_\idot(E) @>{\id+\sigma+\dots+\sigma^{p-1}}>>
  C_\idot(E) @>{\id - \sigma}>> C_\idot(E),\\
@. @AA{B}A @AA{B}A @AA{B}A
\end{CD}
\end{equation}
with $\id$ of degree $(-1,2)$ added to the total differential, where
$C_\idot(E)$, resp. $C'_\idot(E)$ is the complex with terms $E([n])$
and the differential $b_p$, resp. $b'_p$, and $B$ is the horizontal
differential in the complex \eqref{hc.lambda.p} written down for
$p=1$. The complex $C_\idot(E)$ computes the Hochschild homology
$HH_\idot(E)$, the complex $C'_\idot(E)$ is acyclic, and the whole
complex \eqref{hc.hs.eq} computes the cyclic homology $HC_\idot(E)$.

We see that the cyclic homology of the $p$-cyclic object $E$
actually admits two periodicity endomorphisms: $u$ and $u'$. The
horizontal endomorphism $u$ is the usual periodicity map; the
vertical map $u'$ is something new. However, we have the following.

\begin{lemma}\label{no.uprime}
In the situation above, assume that $p = \cchar k$. Then the
vertical periodicity map $u':HC_\idot(E) \to HC_{\idot-2}(E)$ is
equal to $0$.
\end{lemma}

\proof[Sketch of a proof.] It might be possible to write explicitly
a contracting homotopy for the map $u'$, but this is very
complicated; instead, we will sketch the ``scientific'' proof which
uses small categories (for details, see \cite{K2}). For any small
category $\Gamma$, its {\em cohomology} $H^\hdot(\Gamma,k)$ is
defined as
$$
H^\hdot(\Gamma,k) = \Ext^\hdot_{\Fun(\Gamma,k)}(k,k),
$$
where $k$ in the right-hand side is the constant functor. This is an
algebra which obviously acts on $H_\idot(\Gamma,E)$ for any $E \in
\Fun(\Gamma,k)$.

The cohomology $H^\hdot(\Lambda,k)$ of the cyclic category $\Lambda$
is the algebra of polynomials in one generator $u$ of degree $2$, $u
\in H^2(\Lambda,k)$; the action of this $u$ on the cyclic homology
$HC_\idot(-)$ is the periodicity map. The same is true for the
$m$-cyclic categories $\Lambda_m$ for all $m \geq 1$.

Now, recall that we have a natural functor $\pi:\Lambda_p \to
\Lambda$, so that there are two natural elements in
$H^2(\Lambda_p,k)$ -- the generator $u$ and the preimage $\pi^*(u)$
of the generator $u \in H^2(\Lambda,k)$. The action of $u$ gives the
horizontal periodicity endomorphism of the complex \eqref{hc.hs.eq},
and the action of $\pi^*(u)$ gives the vertical periodicity
endomorphism $u'$. We have to prove that if $\cchar k = p$, then
$\pi^*(u) = 0$.

To do this, one uses a version of the Hochschild-Serre spectral
sequence associated to $\pi$ -- namely, we have a spectral sequence
$$
H^\hdot(\Lambda) \otimes H^\hdot(\Z/p\Z,k) \Rightarrow
H^\hdot(\Lambda_p,k).
$$
If $\cchar k = p$, then the group cohomology algebra
$H^\hdot(\Z/p\Z,k)$ is the polynomial algebra $k[u,\eps]$ with two
generators: an even generator $u \in H^2(\Z/p\Z,k)$ and an odd
generator $\eps \in H^1(\Z/p\Z,k)$. Since $H^\hdot(\Lambda_p,k) =
k[u]$, the second differential $d_2$ in the spectral sequence must
send $\eps$ to $\pi^*(u)$, so that indeed, $\pi^*(u)=0$ in
$H^2(\Lambda_p,k)$.
\endproof

Consider now the version of the complex \eqref{hc.hs.eq} which
computes the periodic cyclic homology $HP_\idot(E)$ -- to obtain it,
one has to extend the diagram to the right by periodicity. The rows
of the extended diagam then become the standard complexes which
compute the {\em Tate homology} $\vH_\idot(\Z/p\Z,C_\idot(E))$. We
remind the reader that the Tate homology $\vH_\idot(G,-)$ is a
certain homological functor defined for any finite group $G$ which
combines together homology $H_\idot(G,-)$ and cohomology
$H^\hdot(G,-)$, and that for a cyclic group $\Z/m\Z$ with generator
$\sigma$, the Tate homology $\vH_\idot(\Z/m\Z,W)$ with coefficients
in some representation $W$ may be computed by the $2$-periodic
standard complex
\begin{equation}\label{st.per}
\begin{CD}
@>{d_-}>> W @>{d_+}>> W @>{d_-}>> W @>{d_+}>>
\end{CD}
\end{equation}
with $d_+=\id+\sigma+\dots+\sigma^{m-1}$ and $d_-=\id-\sigma$.

If $W$ is a free module over the group algebra $k[G]$, then the Tate
homology vanishes in all degrees, $\vH_\idot(G,W)=0$. When
$G=\Z/m\Z$, this means that the standard complex is acyclic. If $m$
is prime and equal to the charactertistic of the base field $k$, the
converse is also true -- $\vH_\idot(\Z/m\Z,W)=0$ if and only if $W$
is free over $k[\Z/m\Z]$. We would like to claim a similar vanishing
for Tate homology $\vH_\idot(\Z/p\Z,W_\idot)$ with coefficients in
some {\em complex} $W_\idot$ of $k[\Z/p\Z]$-modules; however, this
is not possible unless we impose some finiteness conditions on
$W_\idot$.

\begin{defn}\label{small.defn}
A complex $W_\idot$ of $k[\Z/p\Z]$-modules is {\em effectively
finite} if it is chain-homotopic to a complex of finite length.  A
$p$-cyclic object $E \in \Fun(\Lambda_p,k)$ is {\em small} if its
standard complex $C_\idot(E)$ is effectively finite.
\end{defn}

Now we can finally state our vanishing result for periodic cyclic
homology.

\begin{prop}\label{no.hp}
Assume that $p = \cchar k$. Assume that a $p$-cyclic object $E \in
\Fun(\Lambda_p,k)$ is small, and that $E([n])$ is a free
$k[\Z/p\Z]$-module for every object $[n] \in \Lambda_p$. Then
$HP_\idot(E)=0$.
\end{prop}

\proof{} To compute $HP_\idot(E)$, let us use the periodic version
of the complex \eqref{hc.hs.eq}. We then have a long exact sequence
of cohomology
$$
\begin{CD}
HP_{\idot-1}(E) @>>> \vH_\idot(\Z/p\Z,C_\idot(E)) @>>> HP_\idot(E)
@>{u'}>>,
\end{CD}
$$
where $\vH_\idot(\Z/p\Z,C_\idot(E))$ is computed by the total
complex of the bicomplex
\begin{equation}\label{Tate.hh}
\begin{CD}
@>{\id+\sigma+\dots+\sigma^{p-1}}>> E([1]) @>{\id-\sigma}>> E([1]) 
@>{\id+\sigma+\dots+\sigma^{p-1}}>> \\
@. @AA{b_p}A @AA{b_p}A\\
@>{\id+\sigma+\dots+\sigma^{p-1}}>> E([2]) @>{\id-\sigma}>> E([2]) 
@>{\id+\sigma+\dots+\sigma^{p-1}}>> \\
@. @AA{b_p}A @AA{b_p}A\\
@. \dots @. \dots @. \\
@. @AA{b_p}A @AA{b_p}A\\
@>{\id+\sigma+\dots+\sigma^{p-1}}>> E([n]) @>{\id-\sigma}>> E([n]) 
@>{\id+\sigma+\dots+\sigma^{p-1}}>> \\
@. @AA{b_p}A @AA{b_p}A.
\end{CD}
\end{equation}
By Lemma~\ref{no.uprime}, the connecting differential in the long
exact sequence vanishes, so that it suffices to prove that
$\vH_\idot(\Z/p\Z,C_\idot(E)) = 0$. Since $E([n])$ is free, all the
rows of the bicomplex \eqref{Tate.hh} are acyclic. But since $E$ is
small, $C_\idot(E)$ is effectively finite; therefore the spectral
sequence of the bicomplex \eqref{Tate.hh} converges, and we are
done.
\endproof

\section{Quasi-Frobenius maps}\label{qfr.sec}

We now fix a perfect base field $k$ of characteristic $p > 0$, and
consider an associative algebra $A$ over $k$. We want to construct a
cyclic-homology version of the Cartier isomorphism \eqref{car.usu}
for $A$. In fact, we will construct a version of the inverse
isomorphism $C^{-1}$; it will be an isomorphism
\begin{equation}\label{car.gen}
C^{-1}:HH_\idot(A)((u))^\tw \longrightarrow HP_\idot(A),
\end{equation}
where, as before, $HH_\idot(A)((u))$ in the left-hand side means
``Laurent power series in one variable $u$ of degree $2$ with
coefficients in $HH_\idot(A)$''.

If $A$ is commutative and $X=\Spec A$ is smooth, then $HH_\idot(A)
\cong \Omega^\hdot(X)$, $HP_\idot(A) \cong H_{DR}^\hdot(X)((u))$,
and \eqref{car.gen} is obtained by inverting \eqref{car.usu} (and
repeating the resulting map infinitely many times, once for every
power of the formal variable $u$). It is known that the commutative
inverse Cartier map is induced by the Frobenius isomorphism; thus to
generalize it to non-commutative algebras, it is natural to start
the story with the Frobenius map.

At first glance, the story thus started ends immediately: the map $a
\mapsto a^p$ is not an algebra endomorphism of $A$ unless $A$ is
commutative (in fact, the map is not even additive, $(x+y)^p \neq
x^p + y^p$ for general non-commuting $x$ and $y$). So, there is no
Frobenius map in the non-commutative world.

However, to analyze the difficulty, let use decompose the usual
Frobenius into two maps:
$$
\begin{CD}
A @>{\phi}>> A^{\otimes p} @>{M}>> A,
\end{CD}
$$
where $\phi$ is given by $\phi(a) = a^{\otimes p}$, and $M$ is the
multiplication map, $M(a_1 \otimes \dots \otimes a_p) = a_1 \cdots
a_p$. The map $\phi$ is very bad (e.g. not additive), but this is
the same both in the commutative and in the general associative
case. It is the map $M$ which creates the problem: it is an algebra
map if and only if $A$ is commutative.

In general, it is not possible to correct $M$ so that it becomes an
algebra map. However, even not being an algebra map, it can be made
to act on Hochschild homology, and we already saw how: we can take
the map \eqref{M.eq} of Subsection~\ref{p.cycl.subs}.

As for the very bad map $\phi$, fortunately, it turns out that it
can be perturbed quite a bit. In fact, the only property of this map
which is essential is the following one.

\begin{lemma}\label{V.otimesp}
Let $V$ be a vector space over $k$, and let the cyclic group
$\Z/p\Z$ act on its $p$-th tensor power $V^{\otimes p}$ by the
cyclic permutation.  Then the map $\phi:V \to V^{\otimes p}$, $v
\mapsto v^{\otimes p}$ sends $V$ into the kernel of either of the
differentials $d_+$, $d_-$ of the standard complex \eqref{st.per}
and induces an isomorphism
$$
V^\tw \to \vH_i(\Z/p\Z,V^{\otimes p})
$$
both for odd and even degrees $i$.
\end{lemma}

\proof{} The map $\phi$ is compatible with the multiplication by
scalars, and its image is $\sigma$-invariant, so that it indeed
sends $V$ into the kernel of either of the differentials
$d_-,d_+:V^{\otimes p} \to V^{\otimes p}$. We claim that it is
additive ``modulo $\Im d_\pm$'', and that it induces an isomorphism
$V^\tw \cong \Ker d_\pm/\Im d_\mp$. Indeed, choose a basis in $V$,
so that $V \cong k[S]$, the $k$-linear span of a set $S$. Then
$V^{\otimes p} = k[S^p]$ decomposes as $k[S^p]=k[S] \oplus k[S^p
\setminus \Delta]$, where $S \cong \Delta \subset S^p$ is the
diagonal. This decomposition is compatible with the differentials
$d_\pm$, which actually vanish on the first summand $k[S]$. The map
$\phi$, accordingly, decomposes as $\phi = \phi_0 \oplus \phi_1$,
$\phi_0:V^\tw \to k[S]$, $\phi_1:V^\tw \to k[S^p \setminus
\Delta]$. The map $\phi_0$ is obviously additive and an isomorphism;
therefore it suffices to prove that the second summand of
\eqref{st.per} is acyclic. Indeed, since the $\Z/p\Z$-action on $S^p
\setminus \Delta$ is free, we have $\vH^\hdot(\Z/p\Z,k[S^p \setminus
\Delta]) = 0$.
\endproof

\begin{defn}\label{quasi.fr.defn}
A {\em quasi-Frobenius map} for an associative unital algebra $A$
over $k$ is a $\Z/p\Z$-equivariant algebra map $F:A^\tw \to
A^{\otimes p}$ which induces the isomorphism
$\vH^\hdot(\Z/p\Z,A^\tw) \to \vH^\hdot(\Z/p\Z,A^{\otimes p})$ of
Lemma~\ref{V.otimesp}.
\end{defn}

Here the $\Z/p\Z$-action on $A^\tw$ is trivial, and the algebra
structure on $A^{\otimes p}$ is the obvious one (all the $p$ factors
commute). We note that since $\vH^i(\Z/p\Z,k) \cong k$ for every
$i$, we have $\vH^i(\Z/p\Z,A^\tw) \cong A^\tw$, so that a
quasi-Frobenius map must be injective. Moreover, since the Tate
homology $\vH^\hdot(\Z/p\Z,A^{\otimes p}/A^\tw)$ vanishes, the
cokernel of a quasi-Frobenius map must be a free $k[\Z/p\Z]$-module.

\medskip

In this Section, we will construct a Cartier isomorphism
\eqref{car.gen} for algebras which admit a quasi-Frobenius map (and
satisfy some additional assumptions). In the interest of full
disclosure, we remark right away that quasi-Frobenius maps are very
rare -- in fact, we know only two examples:
\begin{enumerate}
\item $A$ is the tensor algebra $T^\hdot V$ of a $k$-vector space
  $V$ -- it suffices to give $F$ on the generators, where it
  exists by Lemma~\ref{V.otimesp}.
\item $A=k[G]$ is the group algebra of a (discrete) group $G$ -- a
  quasi-Frobenius map $F$ is induces by the diagonal embedding $G
  \subset G^p$.
\end{enumerate}
However, the general construction of the Cartier map given in
Section~\ref{add.sec} will be essentially the same -- it is only the
notion of a quasi-Frobenius map that we will modify.

\begin{prop}
Assume given an algebra $A$ over $k$ equipped with a quasi-Frobenius
map $F:A^\tw \to A^{\otimes p}$, and assume that the category
$A\bimod$ of $A$-bimodules has finite homological dimension. Then
there exists a canonical isomorphism
$$
\phi:HH_\idot(A)((u)) \cong HP_\idot(A).
$$
\end{prop}

\proof{} Consider the functors $i,\pi:\Lambda_p \to \Lambda$ and the
restrictions 
$$
\pi^*A^\tw_\#,i^*A_\# \in \Fun(\Lambda_p,k).
$$
For any $[n] \in \Lambda_p$, the quasi-Frobenius map $F:A^\tw \to
A^{\otimes p}$ induced a map
$$
F^{\otimes n}:\pi^*A^\tw_\#([n]) = (A^\tw)^{\otimes n} \to
i^*A_\#([n]) = A^{\otimes pn}.
$$
By the definition of a quasi-Frobenius map, these maps commute with
the action of the maps $\tau:[n] \to [n]$ and $m^p_i:[n+1] \to [n]$,
$0 \leq i < n$ (recall that $m_i^p = m_i^{\otimes p}$). Moreover,
since $m^p_0 \circ \tau = m^p_{n+1}$, $F^{\otimes \hdot}$ also
commutes with $m^p_n$. All in all, the collection of the tensor
power maps $F^{\otimes\hdot}$ gives a map $F_\#:\pi^*A_\#^\tw \to
i^*A_\#$ of objects in $\Fun(\Lambda_p,k)$. We denote by $\Phi$
the induced map
$$
\Phi= HP_\idot(F_\#):HP_\idot(\pi^*A^\tw_\#) \to
HP_\idot(\Lambda_p,i^*A_\#).
$$
By Lemma~\ref{i.lemma}, the right-hand side is precisely
$HP_\idot(A)$. As for the left-hand side, we note that $\sigma$ is
trivial on $\pi^*A_\#([n])$ for every $[n] \in \Lambda_p$; therefore
the odd horizontal differentials 
\begin{align*}
\id + \tau + \dots + \tau^{pn-1} &= (\id + \tau + \dots + \tau^{n-1})
\circ (\id + \sigma + \dots + \sigma^{p-1})\\
&= p(\id + \tau + \dots + \tau^{n-1}) = 0
\end{align*}
in \eqref{hc.lambda.p} vanish, and we have
$$
HP_\idot(\pi^*A^\tw_\#) \cong HH_\idot(A^\tw)((u)).
$$
Finally, to show that $\Phi$ is an isomorphism, we recall that the
quasi-Frobenius map $F$ is injective, and its cokernel is a free
$k[\Z/p\Z]$-module. One deduces easily that the same is true for
each tensor power $F^{\otimes n}$; thus $F_\#$ is injective, and its
cokernel $\Coker F_\#$ is such that $\Coker F_\#([n])$ is a free
$k[\Z/p\Z]$-module for any $[n] \in \Lambda_p$. To finish the proof,
use the long exact sequence of cohomology and
Proposition~\ref{no.hp}. The only thing left to check is that
Proposition~\ref{no.hp} is applicable -- namely, that the $p$-cyclic
object $i^*A_\#$ is small in the sense of
Definition~\ref{small.defn}.

To do this, we have to show that the bar complex $C_\idot(A^{\otimes
p},A^{\otimes p}_\sigma)$ which computes $HH_\idot(A^{\otimes
p},A^{\otimes p}_\sigma)$ is effectively finite. It is here that we
need to use the assumption of finite homological dimension on the
category $A\bimod$. Indeed, to compute $HH_\idot(A^{\otimes
p},A^{\otimes p}_\sigma)$, we can choose any projective resolution
$P^p_\idot$ of the diagonal $A^{\otimes p}$-bimodule $A^{\otimes
p}$. In particular, we can take any projective resolution $P_\idot$
of the diagonal $A$-bimodule $A$, and use its $p$-th power. To
obtain the bar complex $C_\idot(A^{\otimes p},A^{\otimes
p}_\sigma)$, one uses the bar resolution $C_\idot'(A)$. However, all
these projective resolutions $P_\idot$ are chain-homotopic to each
other, so that the resulting complexes will be also chain-homotopic
{\em as complexes of $k[\Z/p\Z]$-modules}. By assumption, the
diagonal $A$-bimodule $A$ has a projective resolution $P_\idot$ of
finite length; using it gives a complex of finite length which is
chain-homotopic to $C_\idot(A^{\otimes p},A^{\otimes p}_\sigma)$,
just as required by Definition~\ref{small.defn}.
\endproof

\section{Cartier isomorphism in the general case}\label{add.sec}

\subsection{Additivization.}
We now turn to the general case: we assume given a perfect field $k$
of characteristic $p > 0$ and an associative $k$-algebra $A$, and we
want to construct a Cartier-type isomorphism \eqref{car.gen} without
assuming that $A$ admits a quasi-Frobenius map in the sense of
Definition~\ref{quasi.fr.defn}.

Consider again the non-additive map $\phi:A \to A^{\otimes p}$, $a
\mapsto a^{\otimes p}$, and let us change the domain of its
definition: instead of $A$, let $\phi$ be defined on the $k$-vector
space $k[A]$ spanned by $A$ (where $A$ is considered as a set). Then
$\phi$ obviously uniquely extends to a $k$-linear additive map
\begin{equation}\label{qfr.span}
\phi:k[A] \to A^{\otimes p}.
\end{equation}
Taking the $k$-linear span is a functorial operation: setting $V
\mapsto k[V]$ defines a functor $\Span_k$ from the category of
$k$-vector spaces to itself. The functor $\Span_k$ is non-additive,
but it has a tautological surjective map $\Span_k \to \Id$ onto the
identity functor, and one can show that $\Id$ is the maximal
additive quotient of the functor $\Span_k$. If $V=A$ is an algebra,
then $\Span_k(A)$ is also an algebra, and the tautological map
$\Span_k(A) \to A$ is an algebra map.

We note that in both examples \thetag{i}, \thetag{ii} in
Section~\ref{qfr.sec} where an algebra $A$ did admit a
quasi-Frobenius map, what really happened was that the tautological
surjective algebra map $\Span_k(A) \to A$ admitted a splitting $s:A
\to \Span_k(A)$; the quasi-Frobenius map was obtained by composing
this splitting map $s$ with the canonical map \eqref{qfr.span}.

Unfortunately, in general the projection $\Span_k(A) \to A$ does not
admit a splitting (or at least, it is not clear how to construct
one). In the general case, we will modify both sides of the map
\eqref{qfr.span} so that splittings will become easier to come
by. To do this, we use the general technique of {\em additivization}
of non-additive functors from the category of $k$-vector spaces to
itself.

\medskip

Consider the small category $\V = k\Vect^\fg$ of finite-dimensional
$k$-vector spaces, and consider the category $\Fun(\V,k)$ of {\em
all} functors from $\V$ to the category $k\Vect$ of all $k$-vector
spaces. This is an abelian category. The category $\Fun_{add}(\V,k)$
of all {\em additive} functors from $\V$ to $k\Vect$ is also abelian
(in fact, an additive functor is completely defined by its value at
the one-dimensional vector space $k$, so that $\Fun_{add}(\V,k)$ is
equivalent to the category of modules over $k \otimes_{\Z} k$). We
have the full embedding $\Fun_{add}(\V,k) \subset \Fun(\V,k)$, and
it admits a left-adjoint functor -- in other words, for any functor
$F \in \Fun(\V,k)$ there is an additive functor $F_{add}$ and a map
$F \to F_{add}$ which is universal with respect to maps to additive
functors. This ``universal additive quotient'' is not very
interesting. For instance, if $F$ is the $p$-th tensor power
functor, $V \mapsto V^{\otimes p}$, then its universal additive
quotient is the trivial functor $V \mapsto 0$.

To obtain a useful version of this procedure, we have to consider
the derived category $\D(\V,k)$ of the category $\Fun(\V,k)$ and the
full subcategory $\D_{add}(\V,k) \subset \D(\V,k)$ spanned by
complexes whose homology object lie in $\Fun_{add}(\V,k)$.

The category $\D_{add}(\V,k)$ is closed under taking cones, thus
triangulated (this has to be checked, but this is not difficult),
and it contains the derived category of the abelian category
$\Fun_{add}(\Fun_k,k)$. However, $\D_{add}(\V,k)$ is much larger
than this derived category.\ In fact, even for the identity functor
$\Id \in \Fun_{add}(\V,k) \subset \Fun(\V,k)$, the natural map
$$
\Ext^i_{\Fun_{add}(\V,k)}(\Id,\Id) \to
\Ext^i_{\Fun(\V,k)}(\Id,\Id)
$$
is an isomorphism only in degrees $0$ and $1$. Already in degree
$2$, there appear extension classes which cannot be represented by a
complex of additive functors.

Nevertheless, it turns out that just as for abelian categories, the
full embedding $\D_{add}(\V,k) \subset \D(\V,k)$ admits a
left-adjoint functor. We call it the {\em additivization functor}
and denote by $\Add_\idot:\D(\V,k) \to \D_{add}(\V,k)$.  For any $F
\in \Fun(\V,k)$, $\Add_\idot(F)$ is a complex of functors from $\V$
to $k$ with additive homology functors.

The construction of the additivization $\Add_\idot$ is relatively
technical; we will not reproduce it here and refer the reader to
\cite[Section 3]{K2}. The end result is that first, addivization
exists, and second, it can be represented explicitly, by a very
elegant ``cube construction'' introduced fifty years ago by
Eilenberg and MacLane. Namely, to any functor $F \in \Fun(\V,k)$ one
associates a complex $Q_\idot(F)$ of functors from $\V$ to $k$ such
that the homology of this complex are additive functors, and we have
an explicit map $F \to Q_\idot(F)$ which descends to a universal map
in the derived category $\D(\V,k)$. In fact, the complex
$Q_\idot(F)$ is concentrated in non-negative homological degrees,
and $Q_0(F)$ simply coincides with $F$, so that the universal map is
the tautological embedding $F = Q_0(F) \to Q_\idot(F)$. Moreover,
assume that the functor $F$ is {\em multiplicative} in the following
sense: for any $V,W \in \V$, we have a map
$$
F(V) \otimes F(W) \to F(V \otimes W),
$$
and these maps are functorial and associative in the obvious
sense. Then the complex $Q_\idot(F)$ is also multiplicative. In
particular, if we are given a multiplicative functor $F$ and an
associative algebra $A$, then $F(A)$ is an associative algebra; in
this case, $Q_\idot(F)$ is an associative DG algebra concentrated in
non-negative degrees.

\subsection{Generalized Cartier map.}
Consider now again the canonical map \eqref{qfr.span}. There are two
non-additive functors involved: the $k$-linear span functor $V
\mapsto k[V]$, and the $p$-th tensor power functor $V \mapsto
V^{\otimes p}$. Both are multiplicative. We will denote by
$Q_\idot(V)$ the additivization of the $k$-linear span, and we will
denote by $P_\idot(V)$ the additivization of the $p$-th tensor
power. Since additivization is functorial, the map \eqref{qfr.span}
gives a map
$$
\phi:Q_\idot(V) \to P_\idot(V)
$$
for any finite-dimensional $k$-vector space $V$; if $A=V$ is an
associative algebra, then $Q_\idot(A)$ and $P_\idot(A)$ are
associative DG algebras, and $\phi$ is a DG algebra map. We will
need several small refinements of this construction.
\begin{enumerate}
\item We extend both $Q_\idot$ and $P_\idot$ to arbitrary vector
  spaces and arbitrary algebras by taking the limit over all the
  finite-dimensional subspaces.
\item The $p$-th power $V^{\otimes p}$ carries the permutation
  action of the cyclic group $\Z/p\Z$, and the map \eqref{qfr.span}
  is $\Z/p\Z$-invariant; by the functoriality of the additivization,
  $P_\idot(V)$ also carries an action of $\Z/p\Z$, and the map
  $\phi$ is $\Z/p\Z$-invariant.
\item The map \eqref{qfr.span}, while not additive, respects the
  multiplication by scalars, up to a Frobenius twist; unfortunately,
  the additivization procedure ignores this. From now on, we will
  assume that the perfect field $k$ is actually finite, so that the
  group $k^*$ of scalars is a finite group whose order is coprime to
  $p$. Then $k^*$ acts naturally on $k[V]$, hence also on
  $Q_\idot(V)$, and the map $\phi$ factors through the space
  $\Q_\idot(V) = Q_\idot(V)_{k^*}$ of covariants with respect to
  $k^*$.
\end{enumerate}
The end result: in the case of a general algebra $A$, our
replacement for a quasi-Frobenius map is the canonical map
\begin{equation}\label{gen.qfr}
\phi:\Q_\idot(A)^\tw \to P_\idot(A),
\end{equation}
which is a $\Z/p\Z$-invariant DG algebra map. We can now repeat the
procedure of Section~\ref{qfr.sec} replacing a quasi-Frobenius map
$F$ with this canonical map $\phi$. This gives a canonical map
\begin{equation}\label{gen.car.int}
\Phi:HH_\idot(\Q_\idot(A)_\#)^\tw((u)) \to HP_\idot(P_\idot(A)_\#),
\end{equation}
where $\Q_\idot(A)_\#$ in the left-hand side is a complex of cyclic
objects, and $P_\idot(A)_\#$ in the right-hand side is the complex
of $p$-cyclic objects. There is one choice to be made because both
complexes are infinite; we agree to interpret the total complex
which computes $HP_\idot(E_\idot)$ and $HH_\idot(E_\idot)$ for an
infinite complex $E_\idot$ of cyclic or $p$-cyclic objects as the
{\em sum}, not the {\em product} of the corresponding complexes for
the individual terms $HP_\idot(E_i)$, $HH_\idot(E_i)$.

\medskip

To understand what \eqref{gen.car.int} has to do with the Cartier
map \eqref{car.gen}, we need some information on the structure of DG
algebras $P_\idot(A)$ and $\Q_\idot(A)$.

\medskip

The DG algebra $P_\idot(A)$ has the following structure: $P_0(A)$ is
isomorphic to the $p$-th tensor power $A^{\otimes p}$ of the algebra
$A$, and all the higher terms $P_i(A)$, $i \geq 1$ are of the form
$A^{\otimes p} \otimes W_i$, where $W_i$ is a certain representation
of the cyclic group $\Z/p\Z$. The only thing that will matter to us
is that all the representations $W_i$ are {\em free}
$k[\Z/p\Z]$-modules. Consequently, $P_i(A)$ is free over $k[\Z/p\Z]$
for all $i \geq 1$. For the proofs, we refer the reader to
\cite[Subsection 4.1]{K2}. As a corollary, we see that if $A$ is
such that $A\bimod$ has finite homological dimension, then we can
apply Proposition~\ref{no.hp} to all the higher terms in the complex
$P_\idot(A)_\#$ and deduce that the right-hand of
\eqref{gen.car.int} is actually isomorphic to $HP_\idot(A)$:
$$
HP_\idot(P_\idot(A)_\#) \cong HP_\idot(A).
$$
We note that it is here that it matters how we define the periodic
cyclic homology of an infinite complex (the complex $P_\idot(A)$ is
actually acyclic, so, were we to take the product and not the sum of
individual terms, the result would be $0$, not $HP_\idot(A)$).

The structure of the DG algebra $\Q_\idot$ is more interesting. As
it turns out, the homology $H_i(\Q_\idot(A))$ of this DG algebra in
degree $i$ is isomorphic to $A \otimes \St(k)_i$, where
$\St(k)_\idot$ is the dual to the {\em Steenrod algebra} known in
Algebraic Topology -- more precisely, $\St(k)_i^*$ is the algebra of
stable cohomological operations with coefficients in $k$ of degree
$i$. The proof of this is contained in \cite[Section 3]{K2};
\cite[Subsection 3.1]{K2} contains a semi-informal discussion of why
this should be so, and what is the topological interpretation of all
the constructions in this Section. The topological part of the story
is quite large and well-developed -- among other things, it includes
the notions of {\em Topological Hochschild Homology} and {\em
Topological Cyclic Homology} which have been the focus of much
attention in Algebraic Topology in the last fifteen years. A reader
who really wants to understand what is going on should definitely
consult the sources, some of which are indicated in
\cite{K2}. However, within the scope of the present lectures, we
will leave this subject completely alone. The only topological fact
that we will need is the following description of the Steenrod
algebra in low degrees:
\begin{equation}\label{steen}
\St_i(k)=\begin{cases} k,&\qquad i = 0,1,\\
0,&\qquad 1 < i < 2p-2.
\end{cases}
\end{equation}
The proof can be easily found in any algebraic topology textbook.

Thus in particular, the $0$-th homology of $\Q_\idot(A)$ is
isomorphic to $A$ itself, so that we have an augmentation map
$\Q_\idot(A) \to A$ (this is actually induced by the tautogical map
$Q_0(A) = k[A] \to A$). However, there is also non-trivial homology
in higher degrees. Because of this, the left-hand side of
\eqref{gen.car.int} is larger than the left-hand side of
\eqref{car.gen}, and the canonical map $\Phi$ of \eqref{gen.car.int}
has no chance of being an isomorphism (for a topological
interpretation of the left-hand side of \eqref{gen.car.int}, see
\cite[Subsection 3.1]{K2}).

In order to get an isomorphism \eqref{car.gen}, we have to resort to
splittings again, and it would seem that we gained nothing, since
splitting the projection $\Q_\idot(A) \to A$ is the same as
splitting the projection $\Q_0(A) = k[A]_{k^*} \to A$. Fortunately,
in the world of DG algebras we can get away with something less than
a full splitting map. We note the following obvious fact: any {\em
quasiisomorphism} $f:A_\idot \to B_\idot$ of DG algebras induces an
isomorphism $HH_\idot(A_\idot) \to HH_\idot(B_\idot)$ of their
Hochschild homology. Because of this, it suffices to split the
projection $\Q_\idot(A) \to A$ ``up to a quasiisomorphism''. More
precisely, we introduce the following.

\begin{defn}
A {\em DG splitting} $\langle \overline{A}_\idot,s \rangle$ of a DG
algebra map $f:\wt{A}_\idot \to A_\idot$ is a pair of a DG algebra
$\overline{A}_\idot$ and a DG algebra map $s:\overline{A}_\idot \to
\wt{A}_\idot$ such that the composition $f \circ
s:\overline{A}_\idot \to A_\idot$ is a quasiisomorphism.
\end{defn}

\begin{lemma}\label{car.dgspl}
Assume that the associative algebra $A$ is such that $A\bimod$ has
finite homological dimension. For any DG splitting $\langle
A_\idot,s \rangle$ of the projection $\Q_\idot(A) \to A$, the
composition map
\begin{multline*}
\Phi \circ s:HH_\idot(A)^\tw((u)) \cong HH_\idot(A_\idot)^\tw((u))
\to\\
\to HH_\idot(\Q_\idot(A))^\tw((u)) \to HP_\idot(P_\idot(A)_\#) \cong
HP_\idot(A)
\end{multline*}
is an isomorphism in all degrees.
\end{lemma}

The proof is not completely trivial but very straightforward; we
leave it as an exercise (or see \cite[Subsection 4.1]{K2}). By
virtue of this lemma, all we have to do to construct a Cartier-type
isomorphism \eqref{car.gen} is to find a DG splitting of the
projection $\Q_\idot(A) \to A$.

\subsection{DG splittings.}
To construct DG splittings, we use obstruction theory for DG
algebras, which turns out to be pretty much parallel to the usual
obstruction theory for associative algebras. A skeleton theory
sufficient for our purposes is given in \cite[Subsection
4.3]{K2}. Here are the main points.
\begin{enumerate}
\item Given a DG algebra $A_\idot$ and a DG $A_\idot$-bimodule
  $M_\idot$, one defines {\em Hoch\-schild cohomology}
  $HH^\hdot_{\D}(A_\idot,M_\idot)$ as
$$
HH^\hdot_{\D}(A_\idot,M_\idot) = \Ext^\hdot_{\D}(A_\idot,M_\idot),
$$
where $A_\idot$ in the right-hand side is the diagonal
$A_\idot$-bimodule, and $\Ext^\hdot_{\D}$ are the spaces of maps in
the ``triangulated category of $A_\idot$-bimodules'' -- that is, the
derived category of the abelian category of DG $A_\idot$-bimodules
localized with respect to quasiisomorphisms. Explicitly,
$HH^\hdot(A_\idot,M_\idot)$ can be computed by using the bar
resolution of the diagonal bimodule $A_\idot$. This gives a complex
with terms $\Hom(A_\idot^{\otimes n},M_\idot)$, where $n \geq 0$ is
a non-negative integer, and a certain differential
$\delta:\Hom(A_\idot^{\otimes \hdot},M_\idot) \to
\Hom(A_\idot^{\otimes \hdot+1},M_\idot)$; the groups
$HH^\hdot(A_\idot,M_\idot)$ are computed by the total complex of the
bicomplex
$$
\begin{CD}
M_\idot @>{\delta}>> \Hom(A_\idot,M_\idot) @>{\delta}>> \dots
@>{\delta}>> \Hom(A_\idot^{\otimes \hdots},M_\idot) @>{\delta}>>.
\end{CD}
$$
\item By a {\em square-zero extension} of a DG algebra $A_\idot$ by
  a DG $A_\idot$-bimodule we understand a DG algebra $\wt{A_\idot}$
  equipped with a surjective map $\wt{A_\idot} \to A_\idot$ whose
  kernel is identified with $M_\idot$ (in particular, the induced
  $\wt{A_\idot}$-bimodule structure on the kernel factors through
  the map $\wt{A_\idot} \to A_\idot$). Then square-zero extensions
  are classified up to a quasiisomorphism by elements in the
  Hochschild cohomology group $HH^2_{\D}(A_\idot,M_\idot)$. A
  square-zero extension admits a DG splitting if and only if its
  class in $HH^2_{\D}(A_\idot,M_\idot)$ is trivial.
\end{enumerate}
To apply this machinery to the augmentation map $\Q_\idot(A) \to A$,
we consider the {\em canonical filtration} $\Q_\idot(A)_{\geq
\idot}$ on $\Q_\idot(A)$ defined, as usual, by
$$
\Q_i(A)_{\geq j} =
\begin{cases}
0, &\qquad i \leq j,\\
\Ker d, &\qquad i = j+1,\\
\Q_i(A), &\qquad i > j+1,
\end{cases}
$$
where $d$ is the differential in the complex $\Q_\idot(A)$. We
denote the quotients by $\Q_\idot(A)_{\leq j}(A) =
\Q_\idot(A)/\Q_\idot(A)_{\geq j}$, and we note that for any $j \geq
1$, $\Q_\idot(A)_{\leq j}$ is a square-zero extension of
$\Q_\idot(A)_{\leq j-1}$ by a DG bimodule quasiisomorphic to $A
\otimes \St_j(k)[j]$ (here $\St_j(k)$ is the corresponding term of
the dual Steenrod algebra, and $[j]$ means the degree shift). We use
induction on $j$ and construct a collection $\langle A_\idot^j,s
\rangle$ of compatible DG splittings of the surjections
$\Q_\idot(A)_{\leq j} \to A$. There are three steps.

\smallskip

{\em Step {\normalfont 1}.} For $j=0$, there is nothing to do: the
projection $\Q_\idot(A)_{\leq 0} \to A$ is a quasiisomorphism.

\smallskip

{\em Step {\normalfont 2}.} For $j=1$, it tuns out that the
projection $\Q_\idot(A)_{\leq 1} \to A$ admits a DG splitting if and
only if the $k$-algebra $A$ admits a lifting to a flat algebra
$\wt{A}$ over the ring $W_2(k)$ of second Witt vectors of the field
$k$. In fact, even more is true: DG splittings are in some sense in
a functorial one-to-one correspodence with such liftings; the reader
will find precise statements and explicit detailed proofs in
\cite[Subsection 4.2]{K2}.

\smallskip

{\em Step {\normalfont 3}.} We then proceed by induction. Assume
given a DG splitting $A^j_\idot$, $s:A^j_\idot \to \Q_\idot(A)_{\leq
j}$ of the projection $\Q_\idot(A)_{\leq j} \to A$. Form the ``Baer
sum'' $\overline{A}^j_\idot$ of the map $s$ with the square-zero
extension $p:\Q_\idot(A)_{\leq j+1} \to \Q_\idot(A)_{\leq j}$ --
that is, let
$$
\overline{A}^j_\idot \subset \Q_\idot(A)_{\leq j+1} \oplus A^j_\idot
$$
be the subalgebra obtained as the kernel of the map
$$
\begin{CD}
\Q_\idot(A)_{\leq j+1} \oplus A^j_\idot @>{p \oplus (-s)}>>
\Q_\idot(A)_{\leq j}.
\end{CD}
$$
Then $\overline{A}^j_\idot$ is a square-zero extension of
$A^j_\idot$ by a DG $A^j_\idot$-bimodule $\Ker p$ which is
quasiisomorphic to $A \otimes \St_j(k)[j]$. Since $\A^j_\idot$ is
quasiisomorphic to $A$, these are classified by elements in the
Hochschild cohomology group
$$
HH^2(A^j_\idot,\Ker p) \cong HH^{3+j}(A,A) \otimes \St_{j+1}(k).
$$
If $j < 2p-3$, this group is trivial by \eqref{steen}, so that a
DG splitting $A^{j+1}_\idot$ exists. In higher degrees, we have to
impose conditions on the algebra $A$. Here is the end result.

\begin{prop}\label{DG.spl}
Assume given an associative algebra $A$ over a finite field $k$ of
characteristic $p$ such that
\begin{enumerate}
\item $A$ lifts to a flat algebra over the ring $W_2(k)$ of second
  Witt vectors, and
\item $A\bimod$ has finite homological dimension, and moreover,
  we have $HH^j(A,A) = 0$ whenever $j \geq 2p$.
\end{enumerate}
Then there exists a DG splitting $A_\idot$, $s:A_\idot \to
\Q_\idot(A)$ of the augmentation map $\Q_\idot(A) \to A$.
\end{prop}

\proof{} Construct a compatible system of DG splittings $A^j_\idot$
as described above, and let $A_\idot = \lim_\gets A^j_\idot$.
\endproof

\begin{theorem}\label{car.gen.thm}
Assume given an associative algebra $A$ over a finite field $k$ of
characteristic $p$ which satisfies the assumptions \thetag{i},
\thetag{ii} of Proposition~\ref{DG.spl}. Then there exists an
isomorphism
$$
C^{-1}:HH_\idot(A)((u))^\tw \longrightarrow HP_\idot(A),
$$
as in \eqref{car.gen}.
\end{theorem}

\proof{} Combine Proposition~\ref{DG.spl} and Lemma~\ref{car.dgspl}.
\endproof

This is our generalized Cartier map. We note that the conditions
\thetag{i}, \thetag{ii} that we have to impose on the algebra $A$
are completely parallel to the conditions \thetag{i}, \thetag{ii} on
page \pageref{car.eq} which appear in the commutative case:
\thetag{i} is literally the same, and as for \thetag{ii}, note that
if $A$ is commutative, then the category of $A$-bimodules is
equivalent to the category of quasicoherent sheaves on $X \times X$,
where $X = \Spec A$. By a famous theorem of Serre, this category has
finite homological dimension if and only if $X$ is smooth, and this
dimension is equal to $\dim(X \times X)=2\dim X$.

\section{Applications to Hodge Theory}\label{deg.sec}

To finish the paper, we return to the original problem discussed in
the Introduction: the degeneration of the Hodge-to-de Rham spectral
sequence. On the surface of it, Theorem~\ref{car.gen.thm} is strong
enough so that one can apply the method of Deligne and Illusie in
the non-commutative setting. However, it has one fault. While in the
commutative case we are dealing with an {\em algebraic variety} $X$,
Theorem~\ref{car.gen.thm} is only valid for an associative {\em
algebra}. In particular, were we to try to deduce the classical
Cartier isomorphism \eqref{car.usu} from Theorem~\ref{car.gen.thm},
we would only get it for affine algebraic varieties. In itself, it
might not be completely meaningless. However, the commutative
Hodge-to-de Rham degeneration is only true for a {\em smooth} and
{\em proper} algebraic variety $X$ -- and a variety of dimension
$\geq 1$ cannot be proper and affine at the same time. The general
non-commutative degeneration statement also requires some versions
of properness, and in the affine setting, this reduces to requiring
that the algebra $A$ is finite-dimensional over the base field. A
degeneration statement for such algebras, while not as completely
trivial as its commutative version, is not, nevertheless, very
exciting.

\medskip

Fortunately, the way out of this difficulty has been known for some
time; roughly speaking, one should pass to the level of {\em derived
categories} -- after which all varieties, commutative and
non-commutative, proper or not, become essentially affine.

\medskip

More precisely, one first notices that Hochschild homology of an
associative algebra $A$ is {\em Morita-invariant} -- that is, if $B$
is a different associative algebra such that the category $B\mod$ of
$B$-modules is equivalent to the category $A\mod$ of $A$-modules,
then $HH_\idot(A) \cong HH_\idot(B)$. The same is true for cyclic
and periodic cyclic homology, and for Hochschild cohomology
$HH^\hdot(A)$. In fact, B. Keller has shown in \cite{kel} how to
construct $HC_\idot(A)$ and $HH_\idot(A)$ starting directly from the
abelian category $A\mod$, without using the algebra $A$ at
all.

Moreover, Morita-invariance holds on the level of derived
categories: if there exits a left-exact functor $F:A\mod \to B\mod$
such that its derived functor is an equivalence of the derived
categories $\D(A\mod) \cong \D(B\mod)$, then $HH_\idot(A) \cong
HH_\idot(B)$, and the same is true for $HC_\idot(-)$, $HP_\idot(-)$,
and $HH^\hdot(-)$.

Unfortunately, one cannot recover $HH_\idot(A)$ and other
homological invariants directly from the derived category
$\D(A\mod)$ considered as a triangulated category -- the notion of a
triangulated category is too weak. One has to fix some
``enhancement'' of the triangulated category structure. At present,
it is not clear what is the most convenient choice among several
competing approaches. In practice, however, every ``natural'' way to
construct a triangulated categery $\D$ also allows to equip it with
all possible enhancements, so that the Hochschild homology
$HH_\idot(\D)$ and other homological invariants can be defined.

As long as we work over a fixed field, probably the most convenient
of those ``natural'' ways is provided by the DG algebra
techniques. For every assocative DG algebra $A^\hdot$ over a field
$k$, one defines $HH_\idot(A^\hdot)$, $HC_\idot(A^\hdot)$,
$HP_\idot(A^\hdot)$, and $HH^\hdot(A^\hdot)$ in the obvious way, and
one shows that if two DG algebras $A^\hdot$, $B^\hdot$ have
equivalent triangulated categories $\D(A^\hdot\mod)$,
$\D(B^\hdot\mod)$ of DG modules, then all their homological
invariants such as $HH_\idot(-)$ are isomorphic. Moreover, the DG
algebra approach is versatile enough to cover the case of non-affine
schemes. Namely, one can show that for every quasiprojective variety
$X$ over a field $k$, there exists a DG algebra $A^\hdot$ over $k$
such that $\D(A^\hdot\mod)$ is equivalent to the derived category of
coherent sheaves on $X$. Then $HH_\idot(A^\hdot)$ is the same as the
Hochschild homology of the category of coherent sheaves on $X$, and
the same is true for the other homological invariants -- in
particular, if $X$ is smooth, we have
$$
HH_i(A^\hdot) \cong \bigoplus_j H^j(X,\Omega^{i+j}_X),
$$
and $HC_\idot(A^\hdot)$ is similarly expressed in terms of the de
Rham cohomology groups of $X$. It is in this sense that all the
varieties become affine in the ``derived non-commutative''
world. We note that in general, although $X$ is the usual
commutative algebraic variety, one cannot insure that the algebra
$A^\hdot$ which appears in this construction is also commutative.

\medskip

Thus for our statement on the Hodge-to-de Rham degeneration, we use
the language of associative DG algebras. The formalism we use is
mostly due to B. To\"en; the reader will find a good overview in
\cite[Section 2]{ToVa}, and also in B. Keller's talk \cite{Ke} at
ICM Madrid.

\begin{defn}
Assume given a DG algebra $A^\hdot$ over a field $k$.
\begin{enumerate}
\item $A^\hdot$ is {\em compact} if it is perfect as a complex of
  $k$-vector spaces.
\item $A^\hdot$ is {\em smooth} if it is perfect as the diagonal DG
  bimodule over itself.
\end{enumerate}
\end{defn}

By defintion, a DG $B^\hdot$-module $M_\idot$ over a DG algebra
$B^\hdot$ is perfect if it is a compact object of the triangulated
category $\D(B^\hdot)$ in the sense of category theory -- that is,
we have
$$
Hom(M_\idot,\lim_\to N_\idot) = \lim_\to \Hom(M_\idot,N_\idot)
$$
for any filtered inductive system $N_\idot \in \D(B^\hdot)$. It is
an easy exercise to check that compact objects in the category
$k\Vect$ are precisely the finite-dimensional vector spaces, so that
a complex of $k$-vector spaces is perfect if and only if its
homology is trivial outside of a finite range of degrees, and all
the non-trivial homology groups are finite-dimensional $k$-vector
spaces. In general, there is a theorem which says that a DG module
$M_\idot$ is perfect if and only if it is a retract -- that is, the
image of a projector -- of a DG module $M_\idot'$ which becomes a free
finitely-generated $B^\hdot$-module if we forget the
differential. We refer the reader to \cite{ToVa} for exact
statements and proofs. We note only that if a DG algebra $A^\hdot$
describes an algebraic variety $X$ -- that is, $\D(A^\hdot) \cong
\D(X)$ -- that $A^\hdot$ is compact if and only if $X$ is proper,
and $A^\hdot$ is smooth if and only if $X$ is smooth (for
smoothness, one uses Serre's Theorem mentioned in the end of
Section~\ref{add.sec}).

\begin{theorem}\label{dg.thm}
Assume given an associative DG algebra $A^\hdot$ over a field $K$ of
characteristic $0$. Assume that $A^\hdot$ is smooth and
compact. Moreover, assume that $A^\hdot$ is concentrated in
non-negative degrees. Then the Hodge-to-de Rham spectral sequence
$$
HH_\idot(A^\hdot)[u] \Rightarrow HC_\idot(A^\hdot)
$$
of \eqref{h.dr.eq} degenerates at first term.
\end{theorem}

In this theorem, we have to require that $A^\hdot$ is concentrated
in non-negative degrees. This is unfortunate but inevitable in our
approach to the Cartier map, which in the end boils down to
Lemma~\ref{V.otimesp} -- whose statement is obviously incompatible
with any grading one might wish to put on the vector space $V$. Thus
our construction of the Cartier isomorphism does not work at all for
DG algebras. We circumvent this difficulty by passing from DG
algebras to {\em cosimplicial} algebras -- that is, associative
algebras $\A \in \Fun(\Delta,K)$ in the tensor category
$\Fun(\Delta,K)$ -- for which one can construct the Cartier map
``pointwise'' (it is the passage from DG to cosimplicial algebras
which forces us to require $A^i=0$ for negative $i$). This occupies
the larger part of \cite[Subsection 5.2]{K2}, to which we refer the
reader. Here we will only quote the end result.

\begin{prop}
Assume given a smooth and compact DG algebra $A^\hdot$ over a finite
field $k$ of characteristic $p = \cchar k$. Assume that $A^\hdot$ is
concentrated in non-negative degrees, and that, moreover,
\begin{enumerate}
\item $A^\hdot$ can be lifted to a flat DG algebra over the ring
  $W_2(k)$ of second Witt vectors of the field $k$, and
\item $HH^i(A,A)=0$ when $i \geq 2p$.
\end{enumerate}
Then there exists an isomorphism
$$
C^{-1}:HH_\idot(A^\hdot)((u)) \cong HP_\idot(A^\hdot),
$$
and the Hodge-to-de Rham spectral sequence \eqref{h.dr.eq} for the
DG algebra $A^\hdot$ degenerates at first term.
\end{prop}

As in the commutative case of \cite{DI}, degeneration follows
immediately from the existence of the Cartier isomorphism $C^{-1}$
for dimension reasons. The construction of the map $C^{-1}$
essentially repeats what we did in Section~\ref{add.sec} in the
framework of cosimplicial algebras, with a lot of technical
nuissance because of the need to insure the convergence of various
spectral sequences, see \cite[Subsection 5.3]{K2}. To deduce
Theorem~\ref{dg.thm}, one uses the standard technique of the
reduction to positive characteristic, just as in the commutative
case; this is made possible by the following beautiful theorem due
to B. To\"en \cite{to}.

\begin{theorem}[\cite{to}]
Assume given a smooth and compact DG algebra $A^\hdot$ over a field
$K$. Then there exists a finitely generated subring $R \subset K$
and a DG algebra $A^\hdot_R$, smooth and compact over $R$, such that
$A^\hdot \cong A^\hdot \otimes_R K$.
\end{theorem}

We note that this result does not require the algebra $A^\hdot$ to
be concentrated in non-negative degrees. We expect that neither does
our Theorem~\ref{dg.thm}, but so far, we could not prove it -- the
technical difficulties seem to be much too severe.

\bigskip

\noindent
{\sc
Steklov Math Institute\\
Moscow, USSR}

\bigskip

\noindent
{\em E-mail address\/}: {\tt kaledin@mccme.ru}


\begin{thebibliography}{FPSVW}

\bibitem[Co]{C} A. Connes, {\em Cohomologie cyclique et foncteur
$\Ext^n$}, Comptes Rendues Ac. Sci. Paris S\'er. A-B, {\bf 296}
(1983), 953--958.

\bibitem[DI]{DI} P. Deligne and L. Illusie, {\em Rel\'evements
modulo $p^2$ et d\'ecomposition du complexe de de Rham},
Inv. Math. {\bf 89} (1987), 247--270.

\bibitem[FT1]{FT1} B. Feigin and B. Tsygan, {\em Cohomology of the
Lie algebra of generalized Jacobi matrices}, in Russian,
Funct. An. Appl. {\bf 17} (1983), 86--87.

\bibitem[FT2]{FT} B. Feigin and B. Tsygan, {\em Additive
$K$-Theory}, in Lecture Notes in Math. {\bf 1289} (1987), 97--209.

\bibitem[G2]{Gr} A. Grothendieck, {\em On the de Rham cohomology of
algebraic varieties}, Publ. Math. IHES {\bf 29} (1966), 95-103.

\bibitem[HKR]{HKR} G. Hochschild, B. Kostant, and A. Rosenberg, {\em
Differential forms on regular affine algebras}, Trans. AMS {\bf 102}
(1962), 383--408.

\bibitem[K1]{K1} D. Kaledin, {\em Non-commutative Cartier operator
and Hodge-to-de Rham degeneration}, preprint math.AG/0511665.

\bibitem[K2]{K2} D. Kaledin, {\em Non-commutative Hodge-to-de Rham
degeneration via the method of Deligne-Illusie}, preprint
math.KT/0611623.

\bibitem[K3]{K3} D. Kaledin, {\em Cyclic homology with
coefficients}, preprint\\ math.KT/0702068.

\bibitem[Ke1]{kel} B. Keller, {\em On the cyclic homology of exact
categories}, Journ. of Pure and Appl. Algebra {\bf 136} (1999),
1-56.

\bibitem[Ke2]{Ke} B. Keller, {\em On differential graded categories},
a talk at ICM 2006, Madrid.

\bibitem[KS]{KS} M. Kontsevich and Y. Soibelman, {\em Notes on
A-infinity algebras, A-infinity categories and non-commutative
geometry, I}, preprint math.RA/0606241.

\bibitem[L]{L} J.-L. Loday, {\em Cyclic Homology}, second ed.,
Springer, 1998.

\bibitem[LQ]{LQ} J.-L. Loday and D. Quillen, {\em Homologie cyclique
et homologie de l'alg\`ebre de Lie des matrices}, Comptes Rendues
Ac. Sci. Paris S\'er. A-B, {\bf 296} (1983), 295--297.

\bibitem[T]{to} B. To\"en, {\em Rings of definition of smooth
and proper dg-algebras}, preprint math.AT/0611546.

\bibitem[TV]{ToVa} B. To\"en and M. Vaqui\'e, {\em Moduli of objects
in dg-categories}, preprint math.AG/0503269.

\end{thebibliography}
\end{document}